\documentclass[11pt]{article}
\usepackage{graphicx}
\usepackage{amsfonts}
\usepackage{bbm}
\usepackage{mathrsfs}
 \usepackage{mathrsfs,amsmath,amssymb,amsthm}
\usepackage[dvips]{color}
 \usepackage{amssymb}
 \usepackage{amsbsy}
 \usepackage[dvips]{color}
\allowdisplaybreaks
 \theoremstyle{plain}
\theoremstyle{remark}  \newtheorem{remark}{\noindent\mbox{Remark}}
 \theoremstyle{plain}
 \theoremstyle{plain}
\theoremstyle{plain} \newtheorem{theorem}{\noindent\mbox{Theorem}}
 \theoremstyle{plain}\newtheorem{proposition}{\noindent\mbox{Proposition}}
 \theoremstyle{plain}
\theoremstyle{definition}

 \def\bq{\begin{equation}}
 \def\eq{\end{equation}}
 \def\eqn{\end{eqnarray}}
 \def\bqn{\begin{eqnarray}}
 \def\proof{\noindent{\it Proof.~~}}
 \def\qed{\hfill$\Box$\medskip}
 \def\rto{\rightarrow\infty}
 \def\z{\left}
 \def\y{\right}
 \def\no{\nonumber}

\begin{document}
 \title{\textbf{Birth and death process with one-side bounded jumps in random environment}\footnote{Supported by National
Nature Science Foundation of China (Grant No. 11226199) and Natural Science Foundation of Anhui Educational Committee (Grant No. KJ2014A085)}}                  %%%   the Fund which you are supported by  %%%

\author{  Hua-Ming \uppercase{Wang} \\
     Department of Mathematics, Anhui Normal University, Wuhu 241003, China\\
    E-mail\,$:$ hmking@mail.ahnu.edu.cn}
\date{}
\maketitle%

\vspace{-.7cm}

\begin{center}
\begin{minipage}[c]{12cm}
\begin{center}\textbf{Abstract}\quad \end{center}

Let $\omega=(\omega_i)_{i\in\mathbb Z}=(\mu^{L}_i,...,\mu^{1}_i,\lambda_i)_{i\in \mathbb Z}$, which serves as the environment, be a sequence of i.i.d. random nonnegative vectors, with  $L\ge1$  a positive integer. We study birth and death process $N_t$ which, given the environment $\omega,$ waits at a state $n$ an exponentially distributed time with parameter $\lambda_n+\sum_{l=1}^L\mu^{l}_n$ and then jumps to $n-i$ with probability ${\mu^i_n}/(\lambda_n+\sum_{l=1}^L\mu^{l}_n),$ $i=1,...,L$ or to $n+1$ with probability ${\lambda_n}/(\lambda_n+\sum_{l=1}^L\mu^{l}_n).$ A sufficient condition for the existence,  a criterion for  recurrence, and  a law of large numbers of the process $N_t$ are presented.  We show that the first passage time $T_1\overset{\mathscr D}{=}\xi_{0,1}+\sum_{i\le -1}\sum_{k=1}^{U_{i,1}}\xi_{i,k}+\sum_{i\le -1}\sum_{k= 1}^{U_{i,1}+...+U_{i,L}}\tilde{\xi}_{i+1,k},$ where $(U_{i,1},...,U_{i,L})_{i\le0}$ is an $L$-type branching process in random environment and, given $\omega,$ $\xi_{i,k},\ \tilde\xi_{i,k},\ i\le 0,\ k\ge 1$ are mutually independent random variables such that
  $P_\omega(\xi_{i,k}\ge t)=e^{-(\lambda_i+\sum_{l=1}^L\mu^{l}_i)t},\ t\ge 0.$ This fact enables us to give an explicit velocity of the law of large numbers.

\mbox{}\textbf{Keywords:}\  birth and death process; random environment; the first passage time; branching structure.\\
\mbox{}\textbf{MSC 2010:}\
 60K37;  60J80
\end{minipage}
\end{center}

\section{Introduction}

\subsection{Model and background}
The aim of this paper is to study the birth and death process with one-side bounded jumps in random environment. To construct the environment, fix $1\le L\in \mathbb Z$ and let $\Omega$ be the collection of $\omega=(\omega_i)_{i\in\mathbb Z}=(\mu^{L}_i,...,\mu^{1}_i,\lambda_i)_{i\in \mathbb Z},$ where $\lambda_i,\mu_i^l\ge 0$ for all $i\in\mathbb Z,\ l=1,..,L.$ Equip $\Omega$ with the Borel $\sigma$-algebra $\mathcal{F}$ and let $\mathbb P$ be a probability measure on $(\Omega,\mathcal F)$ which makes $(\omega_i)_{i\in\mathbb Z}$ a sequence of i.i.d. random vectors. Then the so-called  random environment is a random element of $\Omega$ chosen according to $\mathbb P.$
Given a realization of $\omega,$ let $N_t$ be a continuous time Markov chain, which waits at a state $n$ an exponentially distributed time with parameter $\lambda_n+\sum_{l=1}^L\mu^{l}_n$ and then jumps to $n-i$ with probability ${\mu^i_n}/(\lambda_n+\sum_{l=1}^L\mu^{l}_n),$ $i=1,...,L$ or to $n+1$ with probability ${\lambda_n}/(\lambda_n+\sum_{l=1}^L\mu^{l}_n).$  We call the process $N_t$ a {\it Birth and Death Process in Random Environment} (BDPRE hereafter) with bounded jumps.

Such a process is the continuous time analogue of a special case of {\it Random Walk  in Random Environment} (RWRE hereafter) with bounded jumps which was introduced in Key \cite{key} and further developed in Letchikov \cite{letb,letc}, Bremont \cite{bra, brc}, Hong and Zhang \cite{hz},  Hong and Wang \cite{hwa,hwb} etc.

The nearest neighbour setting ($L=1$) was studied in Ritter \cite{rit}, where the existence, the criteria for recurrence, and the law of large numbers (LLN hereafter) of the process were studied. The work of Ritter \cite{rit} could be carried out because the birth and death process (with jump size exactly one) was well developed. However, for birth and death process with bounded jumps, in the literatures we are aware of, few result was known. Therefore, to study BDPRE with bounded jumps, one needs to study the behaviors of  birth and death process with bounded jumps.

In this paper, by a classical argument of the existence and uniqueness of the $Q$-process, we give a sufficient condition which ensures the existence of BDPRE with bounded jumps. Then  criterion for recurrence of the process $N_t$ is presented, which depends on the counterpart of RWRE with bounded jumps. Finally, we prove the LLN of $N_t.$

In order to figure out the explicit asymptotic velocity of LLN, we study the first passage time $T_1:=\inf[t>0:N_t>0]$  of  $N_t.$ The idea is as follows. By looking at only the discontinuities of $N_t,$ we get its embedded process $\chi_n,$ which is a discrete time RWRE with bounded jumps. By the branching structure of $\chi_n$ derived in \cite{hwa}, one could use a multitype branching process in random environment to count exactly how many times $N_t$ has ever visited state $i$ before $T_1.$ But every time it visits $i,$ it would wait here an exponentially distributed time period. In this way, we could decompose $T_1$ and consequently give the explicit velocity for the LLN.

\subsection{Main results}

 For a typical realization of $\omega,$ $P_{\omega}$ denotes the law induced by the process $N_t$ starting from $0.$ The measure $P_{\omega}$ is usually related as the {\it quenched} probability. The so-called {\it annealed} probability $P$ is defined by $P(\cdot)=\int_{\Omega}P_{\omega}(\cdot)\mathbb P(d\omega).$
The notations $E_{\omega},$ $E$ and $\mathbb E$ will be used to denote the expectation operators with respect to $P_{\omega},$ $P$ and $\mathbb P$ respectively.

Set for $i,j\in\mathbb Z,$ $$q_{ij}=\left\{\begin{array}{ll}
                                       \lambda_i,&\text{ if } j=i+1; \\
                                       \mu_{i}^l,&\text{ if } j=i-l,\ l=1,...,L;\\
                                       -\big(\lambda_i+\sum_{l=1}^L\mu_i^l), &\text{ if } j=i;\\
                                       0, &\text{ else,}
                                     \end{array}
\right.$$
and let $Q=(q_{ij})$  which is obviously a conservative $Q$-matrix.

Given $\omega,$ consider the conservative $Q$-matrix $Q=(q_{ij}).$
One follows from classical argument that there exists at least one transition matrix $(\overline{p}_{ij}(t))$ such that
\begin{equation}\label{uniq}
  \lim_{t\rightarrow0}\frac{\overline p_{ij}(t)-\delta_{ij}}{t}=q_{ij},\ i,j\in\mathbb Z.
\end{equation}
Net we give a sufficient condition to ensure  such $(\overline{p}_{ij}(t))$ to be unique.

\noindent  {\bf(C1)} $\mathbb P\left(\lambda_0+\sum_{l=1}^L\mu_0^l>0\right)=1;$

\noindent \textbf{(C2)} $\mathbb P\Big(\displaystyle\sum_{n=1}^\infty\frac{1}{\lambda_n+\sum_{n=1}^L\mu_n^l}=\infty,\ \sum_{n=-\infty}^0\frac{1}{\max_{1\le k\le L}\{\lambda_{nL-k}+\sum_{l=1}^L\mu_{nL-k}^l\}}=\infty\Big)=1.$

\begin{proposition}[Existence of $N_t$]\label{exsit}
  Suppose that conditions (C1) and (C2) are satisfied.
Then,  $\mathbb P$-a.s.,  there is only one transition matrix $(\overline{p}_{ij}(t))$ which solves (\ref{uniq}).
\end{proposition}
Proposition \ref{exsit} says that under conditions (C1) and (C2), $\mathbb P$-a.s., the $Q$-process $N_t$ exists. Next we give criteria for the recurrence and transience of $N_t.$
Introduce matrices $$M_i=\left(
    \begin{array}{cccc}
   \frac{\mu_i^1}{\lambda_i} & ... & \frac{\mu_i^{L-1}}{\lambda_i} & \frac{\mu_i^L}{\lambda_i} \\
   1+\frac{\mu_i^1}{\lambda_i} & ... & \frac{\mu_i^{L-1}}{\lambda_i} & \frac{\mu_i^L}{\lambda_i} \\
    \vdots & \ddots & \vdots & \vdots \\
    \frac{\mu_i^1}{\lambda_i} & ... & 1+\frac{\mu_i^{L-1}}{\lambda_i} & \frac{\mu_i^L}{\lambda_i}\\
    \end{array}
     \right),\ i\in\mathbb Z.
 $$
To study the asymptotic behaviour of $M_0M_1\cdots M_n$ we need condition

\noindent\textbf{(C3)} $\mathbb E \ln \frac{\lambda_0}{\lambda_0+\sum_{l=1}^L\mu^{l}_0}>-\infty, \mathbb E \ln \frac{\mu_0^L}{\lambda_0+\sum_{l=1}^L\mu^{l}_0}>-\infty. $

Since $M_i$ depends only on $\omega_i,$ $(M_i)_{i\in \mathbb Z}$ is a sequence of i.i.d. random matrices under $\mathbb P.$ Under condition (C3), $\mathbb E |\ln \|M_0^{-1}\||+\mathbb E |\ln \|M_0\||<\infty.$ Hence one could use Oseledec's multiplicative ergodic theorem (see \cite{osel}) to the sequence $(M_i)_{i\in \mathbb Z}.$ Consequently, we get the Lyapunov
exponents of the sequence $(M_i)_{i\in\mathbb Z}$ which we write in increasing order as
$$-\infty<\gamma_1\le \gamma_2\le...\le\gamma_L<\infty.$$
\begin{theorem}[Recurrence criteria]\label{rtcrit} Suppose that conditions (C1-C3) are all satisfied. Let $\gamma_1\le \gamma_2\le...\le\gamma_{L}$ be the Lyapunov exponents of the sequence $(M_i)_{i\in\mathbb Z}.$ Then

\noindent $\gamma_L< 0\Rightarrow P(\lim_{t\rto}N_t=\infty)=1;$

\noindent $\gamma_L=0\Rightarrow  P(-\infty=\liminf_{t\rto}N_t<\limsup_{t\rto}N_t=\infty)=1;$

\noindent $\gamma_L>0\Rightarrow  P(\lim_{t\rto}N_t=-\infty)=1.$
\end{theorem}

Next we study the asymptotic velocity of the process $N_t.$ Let $T_0=0$ and  define recursively $$T_{n}=\inf\{t>0:N_t=n\}$$ for $n\ge1.$
$T_n$ is the first passage time of $n$ by the process $N_t.$

By Theorem \ref{rtcrit}, if $\gamma_L\le 0,$ $N_t$ is either transient to the right or recurrent. We have
\begin{theorem}[Decomposition of $T_1$]\label{fp}
  Suppose that conditions (C1-C3) are all satisfied and $\gamma_L\le 0.$ Then $P(T_1<\infty)=1$ and
$$T_1\overset{\mathscr D}{=}\xi_{0,1}+\sum_{i\le -1}\sum_{k=1}^{U_{i,1}}\xi_{i,k}+\sum_{i\le -1}\sum_{k= 1}^{U_{i,1}+...+U_{i,L}}\tilde{\xi}_{i+1,k},$$
  where $``\overset{\mathscr D}{=}"$ means ``equal in distribution", $(U_i)_{i\le 0}$ is an $L$-type branching process in random environment whose offspring distributions are given as (\ref{up01}) and (\ref{up02}) below, and given $\omega,$ $ \xi_{i,k},\ \tilde\xi_{i,k},\ i\le 0,\ k\ge 1$ are mutually independent random variables, which are also all independent of $(U_i)_{i\le0},$ such that   $P(\xi_{i,k}\ge t)=e^{-\z(\lambda_i+\sum_{i=1}^L\mu_i^l\y)t},\ t\ge 0.$
  Moreover, with empty product being identity, \begin{equation}\label{et}
    E_\omega T_1=\sum_{i=-\infty}^0\frac{1}{\lambda_i}\mathbf{e}_1M_0M_{-1}\cdots M_{i+1}\mathbf1
  \end{equation} where $\mathbf1=(\mathbf{e}_1+\mathbf{e}_2+...+\mathbf{e}_L)^T,$ and for $1\le i\le L,$ $\mathbf{e}_i$ is a row vector with the $i$th component $1$ and all other components $0.$
\end{theorem}
\begin{remark}
To proof Theorem \ref{fp}, the idea is as follows. By the branching structure for (L,1) RWRE set up in \cite{hwa}, one could count exactly how many times $N_t$ has ever visited $i$ before $T_1.$ Every time it visits $i,$ it will wait here for an exponentially distributed time period. By this approach, we could decompose $T_1$ and study its distribution. $(ET_1)^{-1}$ serves as the asymptotic velocity of $N_t.$
\end{remark}
Define $S(\omega):=\frac{1}{\lambda_0}\sum_{n=1}^\infty \mathbf{e}_1M_1M_{2}\cdots M_{n}\mathbf1.$
\begin{theorem}[LLN of $N_t$]\label{lln}
  Suppose that conditions (C1-C3) are all satisfied and $\gamma_L\le 0.$  Then

\noindent(a) $\mathbb ES(\omega)<\infty\Rightarrow \lim_{t\rto}\frac{N_t}{t}=\mathbb (ES(\omega))^{-1},$ $P$-a.s.;

\noindent(b) $\mathbb ES(\omega)=\infty\Rightarrow \lim_{t\rto}\frac{N_t}{t}=0,$ $P$-a.s..
\end{theorem}

The left part of the paper is arranged as follows. In Section \ref{pec}, we study the existence of the process $N_t$ and give its recurrence/transience criteria. Then in Section \ref{plt}, we give the proof of the LLN and study the distribution of the ladder time $T_1.$

\section{The existence and the recurrence criteria of $N_t$ }\label{pec}

{\it2.1 The existence-Proof of Proposition \ref{exsit}}

 Given $\omega,$ from the classical argument of continuous time Markov chain, $Q$-transitional probability matrix always exits.  Let $(\overline{p}_{ij}(t))$ be a transition matrix which solves (\ref{uniq}) and let $(Y_t)_{t\ge0}$ be the $Q$-process with transition matrix $(\overline{p}_{ij}(t)).$ Let $\tau_0=0$ and define recursively for $n\ge 1,$ $\tau_n=\inf\{t\ge\tau_{n-1}:Y_t\neq Y_{\tau_{n-1}}\}$ where we use the convention $\inf\phi=\infty.$ Then $\tau_n,n=1,2,...$ are consecutive time of discontinuities of $(Y_t)_{t\ge 0}.$

Conditioned on $\{\tau_{n-1}<\infty\}$ and $\{Y_{\tau_{n-1}}=j\},$  if $q_j:=-q_{jj}=\lambda_j+\sum_{l=1}^L\mu_{j}^l>0,$ then $\tau_n-\tau_{n-1}$ has exponential distribution with parameter $q_{j}.$ Therefore we have that $\tau_n<\infty.$ Consequently (C1) implies that for all $n,$ $P$-a.s., $$0=\tau_1<\tau_2<...<\tau_n<\infty.$$
Let $\chi_n=Y_{\tau_n}.$ Then $(\chi_n)_{n\ge 0}$ forms a discrete time Markov chain with transition matrix $(r_{ij})$ whose entries
\begin{equation}\label{rij}
  r_{ij}=\left\{\begin{array}{cl}
                  \frac{\lambda_i}{\lambda_i+\sum_{l=1}^L\mu_i^l}, & j=i+1 \\
                 \frac{\mu_i^l}{\lambda_i+\sum_{l=1}^L\mu_i^l}, & j=i-l, l=1,..,L\\
                 0,& \text{else.}
                \end{array}
\right.
\end{equation}
If \begin{equation}\label{tin}
  P\Big( \sum_{n=0}^\infty q^{-1}_{\chi_n}=\infty\Big)=1,
\end{equation} then we have (see Chung \cite{chung}, Theorem 1 in II.19) that $P(\lim_{n\rto}\tau_n=\infty)=1,$ which implies  the $\mathbb P$-a.s. uniqueness of  $(\overline{p}_{ij}(t)).$
Next we show that  (C2) implies (\ref{tin}). In fact, if the process $(\chi_n)_{n\ge 0}$ is recurrent or transient to the right, then $P$-a.s.,
$$\sum_{n=0}^\infty q^{-1}_{\chi_n}\ge \sum_{n=1}^\infty\frac{1}{\lambda_n+\sum_{l=1}^L\mu_n^l}=\infty.$$
Else if the process $(\chi_n)_{n\ge 0}$ is transient to the left, it must visit at least one state of each of the sets $A_n:=\{nL-k\}_{k=1}^L,n=0,-1,-2,....$  It follows that $ P$-a.s., $$\sum_{n=0}^\infty q^{-1}_{\chi_n}\ge \sum_{n=-\infty}^0\frac{1}{\max_{1\le k\le L}\{\lambda_{nL-k}+\sum_{l=1}^L\mu_{nL-k}^l\}}=\infty.$$ Consequently (\ref{tin}) follows.\qed

\noindent{\it 2.2 Recurrence criteria-Proof of Theorem \ref{rtcrit}}

Under conditions (C1) and (C2), it follows from Proposition \ref{exsit} that the BDPRE with bounded jumps $N_t$ exists.
Let $\tau_0=0,\ \tau_{n}, n\ge 1$ be the consecutive discontinuities of $N_t.$ Set $\chi_n=N_{\tau_n}.$ Given $\omega,$ $\chi_n$ is a discrete time random walk  with transition probabilities $r_{ij}$ defined in (\ref{rij}). $\chi_n$ is also known as the embedded process of $N_t.$
Note that under probability $P,$ $N_t$ and $\chi_n$ have the same recurrence criteria.
Thus Theorem \ref{rtcrit} follows from the following theorem which is a corollary of Theorem A in Letchikov \cite{letc}.
\begin{theorem}\label{rsk} Suppose that conditions (C1-C3) are all satisfied. Let $\gamma_1\le \gamma_2\le...\le\gamma_{L}$ be the Lyapunov exponents of the sequence $(M_i)_{i\in\mathbb Z}.$ Then

\noindent $\gamma_L< 0\Rightarrow P(\lim_{n\rto}\chi_n=\infty)=1;$

\noindent $\gamma_L=0\Rightarrow  P(-\infty=\liminf_{n\rto}\chi_n<\limsup_{t\rto}\chi_n=\infty)=1;$

\noindent $\gamma_L>0\Rightarrow  P(\lim_{n\rto}\chi_n=-\infty)=1.$
\end{theorem}
\proof For $i\in \mathbb Z,$ let $a_i(k)=\frac{\sum_{l=k}^L\mu_i^l}{\lambda_i},$ $k=1,...,L,$  $b_i(1)=\frac{\lambda_i}{\mu_i^L}$ and $b_i(k)=\frac{\sum_{l=k-1}^{L}\mu_{i}^l}{\mu_i^L},$ $k=2,...,L.$
Introduce matrices
\begin{equation}B_i=\left(
        \begin{array}{cccc}
          0 & 1 & \cdots & 0 \\
           \vdots& \vdots & \ddots &\vdots  \\
          0 & 0 & \cdots & 1 \\
          b_{i}(1) & -b_i(2) & \cdots & -b_{i}(L) \\
        \end{array}
      \right) \text{ with } B_i^{-1}=\left(
        \begin{array}{cccc}
          a_{i}(1) &\cdots  & a_i(L-1)& a_i(L) \\
          1 & \cdots & 0 & 0 \\
          \vdots& \ddots & \vdots &\vdots  \\
          0 & \cdots & 1 & 0 \\
        \end{array}
      \right).
\end{equation}
Since $B_i$ depends only on $\omega_i,$ $(B_{i})_{i\in\mathbb Z}$ is a sequence of i.i.d. random matrices under $\mathbb P.$  Under condition (C3) we have that $$\mathbb E|\log \|B_i\||<\infty,\ \mathbb E|\log \|B_i^{-1}\||<\infty.$$
Therefore we can use Oseledec's multiplicative ergodic theorem to get the Lyapunov exponents $(B_i)_{i\in\mathbb Z}$ which we write in increasing order as
$$-\infty<\zeta_1(B)\le \zeta_2(B)\le...\le\zeta_L(B)<\infty.$$
 And those Lyapunov exponents for $(B_i^{-1})_{i\in\mathbb Z}$ are $$-\infty<-\zeta_L(B)\le -\zeta_{L-1}(B)...\le -\zeta_1(B)<\infty.$$
In Theorem A of Letchikov \cite{letc}, the author showed that $P$-a.s., $\chi_n$ is transient to the right, recurrent or transient to the left according as $\zeta_1(B)>0,$ $\zeta_1(B)=0$ or $\zeta_1(B)<0.$ Therefore, if we can show that $\gamma_L=-\zeta_1(B),$ then Theorem \ref{rsk} follows. Indeed, since, for $n\ge L,$ all entries of the product $B_1^{-1}B_2^{-1}\cdots B_{n}^{-1}$ are strictly positive, we have that, as the top Lyapunov exponent of $(B_i^{-1})_{i\in\mathbb Z},$
\begin{equation}\label{topb}
  -\zeta_1(B)=\lim_{n\rto}\mathbb E\log\|B_1^{-1}\cdots B_{n}^{-1}\|.
\end{equation}
Let $$\Lambda=\left(
                \begin{array}{cccc}
                  1 &  &  &  \\
                  1 & 1 &  &  \\
                  \vdots & \vdots & \ddots &  \\
                  1 & 1 & \cdots & 1 \\
                \end{array}
              \right)\text{ with } \Lambda^{-1}=\left(
                \begin{array}{cccc}
                  1 &  &  &  \\
                  -1 & 1 &  &  \\
                  & \ddots & \ddots &  \\
                   &  & -1 & 1 \\
                \end{array}
              \right).
$$
Then we have that $$B_1^{-1}\cdots B_n^{-1}=\Lambda^{-1}M_1\cdots M_n\Lambda.$$
Substituting to (\ref{topb}), it follows that,
$$-\zeta_1(B)=\lim_{n\rto}\mathbb E\log\|\Lambda^{-1}M_1\cdots M_n\Lambda\|=\lim_{n\rto}\mathbb E\log\|M_1\cdots M_n\|=\gamma_L,$$
where the last equality holds because all entries of the products $M_1\cdots M_n$ are strictly positive.
Then Theorem \ref{rsk} is proved.\qed

\section{LLN and the first passage time}\label{plt}

\noindent{\it 3.1 The first passage time-Proof of Theorem \ref{fp}:}

Suppose that conditions (C1) and (C2) hold. Then for $\mathbb P$-a.a. $\omega$  there exists a unique standard transition matrix $(p_{ij}(t))$ which solves (\ref{uniq}).
Let $(N_t)_{t\ge 0}$ be a continuous time Markov chain with standard transition matrices $(p_{ij}(t)).$ Then
\begin{equation}
  \begin{split}
    &P_\omega(N_{t+h}=i+1|N_t=i)=\lambda_ih+o(h);\\
    &P_\omega(N_{t+h}=i-l|N_t=i)=\mu_lh+o(h),\ l=1,2,...,L;\\
    &P_\omega(N_{t+h}=i|N_t=i)=1-\Big(\lambda_ih+\sum_{l=1}^L\mu_lh\Big)+o(h).\\
  \end{split}
\end{equation}
Let $\tau_0=0,$ and $\tau_n=\inf\{t\ge\tau_{n-1}:N_t\neq \tau_{n-1}\}$ for $n\ge1.$ Set $\chi_n:=N_{\tau_n}.$ Then $(\chi_n)_{n\ge0}$ forms a discrete time Markov chain with transition matrix $(r_{ij})$ defined in (\ref{rij}).

%%%%%%%%%%%%%%%%%%%%%%%%%%%%%%%%%%%%%%%%%%%%%%Preliminary result%%%%%%%%%%%%%%%%%%%%%%%%%%%%%%%%%%%%%%%%%%

 For $n\ge 0,$ define $T_n=\inf\{t\ge 0:N_t=n\},$ being the first passage time of state $n$ by $N_t.$
 Next we study  the distribution and the mean of $T_1.$

If $\gamma_L\le 0,$ then by Theorem \ref{rtcrit}, $P$-a.s., both $N_t$ and $\chi_n$ are either recurrent or transient to the right. One follows that $P(T_1<\infty)=1.$ Considering $\chi_n,$ let $\overline{T}_1=\inf\{k>0:\chi_k=1\}.$ Then  $P(\overline T_1<\infty)=1.$ Set $U_0=\mathbf{e}_1,$ and define, for $-\infty<i<0,$ $1\le l\le L,$
 $$U_{i,l}=\#\{0<k<\overline T_1:\chi_{k-1}>i,\chi_k=i-l+1\}$$ and
set
$$U_i:=(U_{i,1},U_{i,2},\cdots,U_{i,L}).$$ Here and throughout, $``\#\{\ \}"$ denotes the number of elements in set $``\{\ \}".$
Note that $U_{i,1}$ is the total number of steps by $\chi_n$ which jumps downwards from some state above $i$ to $i$ before $\overline{T}_1$ and $U_{i,2}+...+U_{i,L}$ is the total number of steps by $\chi_n$ which cross $i$ downwards before $\overline{T}_1.$ Since $\tau_n,n=0,1,2,...$ are consecutive discontinuities of $N_t,$ the total number of negative jumps  of $N_t$ which reach $i$ equals to $U_{i,1}$ and the total number of negative jumps of $N_t$ which cross $i$ downwards before $T_1$ equals to $U_{i,2}+...+U_{i,L}.$

Suppose that a particle moves along the path of $N_t.$ Firstly, the particle starts from $0$ and it stays at $0$ for a time period $\xi_{0,1}$ with $P(\xi_{0,1}\ge t)=e^{-\big(\lambda_0+\sum_{l=1}^L\mu_0^l\big)t},\ t\ge 0.$

Secondly we consider the waiting time caused by the negative jumps before $T_1$.
After the $k(\ge1)$-th visit of $i$ by a negative jump, the particle will stay at $i$ with a random time $\xi_{i,k}$ with $P(\xi_{i,k}\ge t)=e^{-\big(\lambda_i+\sum_{l=1}^L\mu_i^l\big)t},\ t\ge 0.$ Then the total amount of time the particle stays at $i$ caused by those negative jumps which reach $i$ downwards before $T_1$ is $\sum_{k=1}^{U_{i,1}}\xi_{i,k}.$ The total amount of time that the particle stays at the negative half lattice caused by those negative jumps before $T_1$ equals to $\sum_{i\le -1}\sum_{k=1}^{U_{i,1}}\xi_{i,k}.$  By the strong Markov property, $\xi_{i,k},\ i\le 0,\ k\ge 1$ are mutually independent and $U_i, \xi_{i,k},k\ge 1$ are also mutually independent.

Thirdly, we consider the waiting time caused by the positive jumps before $T_1.$ Since $N_{T_1}=1,$ and the positive jumps are nearest neighbor, once the particle takes a negative jump downwards from some state above $i$ to $i$ or across $i,$ it has to take a positive jump from $i$ to $i+1$ in order to reach the state $1$ finally.  In this point of view, we have that the number of jumps of the particle before $T_1$ from $i$ to $i+1$ equals to $U_{i,1}+...+U_{i,L}.$
After the $k(\ge1)$-th visit of $i+1$ by a positive jump, the particle will stay at $i+1$ with a random time $\tilde\xi_{i+1,k}$ with $P(\tilde\xi_{i+1,k}\ge t)=e^{-\big(\lambda_{i+1}+\sum_{l=1}^L\mu_{i+1}^l\big)t},\ t\ge 0.$ Then the total amount of time the particle stays at $i+1$ caused by those jumps from $i$ to $i+1$  before $T_1$ is $\sum_{k=1}^{U_{i,1}+...+U_{i,L}}\tilde\xi_{i+1,k}.$ The total amount of time the particle stays at the negative half lattice caused by those positive jumps before $T_1$ equals to $\sum_{i\le -1}\sum_{k=1}^{U_{i,1}+...+U_{i,L}}\tilde\xi_{i+1,k}.$  By the strong Markov property, $\xi_{i,k},\ \tilde\xi_{i,k},\ i\le 0,\ k\ge 1$ are mutually independent and they are all independent of $U_i.$

The above discussion yields that
$$T_1\overset{\mathscr D}{=}\xi_{0,1}+\sum_{i\le -1}\sum_{k=1}^{U_{i,1}}\xi_{i,k}+\sum_{i\le -1}\sum_{i= 1}^{U_{i,1}+...+U_{i,L}}\tilde{\xi}_{i+1,k}.$$

 On the other hand, in Hong and Wang \cite{hwa}, Theorem 1.1, the authors showed that $(U_n)_{n\le0}$ forms a multitype branching process with offspring distribution \begin{eqnarray}\label{up01}
&&P_\omega(U_{i-1}=(u_1,...,u_L)\big|U_{i}=\mathbf{e}_1)\no\\
&&\quad\quad=\frac{(u_1+...+u_L)!}{u_1!\cdots
u_L!}\z(\frac{\mu_i^1}{\lambda_i+\sum_{l=1}^L\mu_i^l}\y)^{u_1}\cdots\z(\frac{\mu_i^L}{\lambda_i+\sum_{l=1}^L\mu_i^l}\y)^{u_L}\z(\frac{\lambda_i}{\lambda_i+\sum_{l=1}^L\mu_i^l}\y),
\end{eqnarray}
and for $2\le l\le L,$
\begin{eqnarray}\label{up02}
&&P_\omega\z(U_{i-1}=(u_1,...,1+u_{l-1},...,u_L)\big|U_{i}=\mathbf{e}_l\y)\no\\
 &&\quad\quad=\frac{(u_1+...+u_L)!}{u_1!\cdots
u_L!}\z(\frac{\mu_i^1}{\lambda_i+\sum_{l=1}^L\mu_i^l}\y)^{u_1}\cdots\z(\frac{\mu_i^L}{\lambda_i+\sum_{l=1}^L\mu_i^l}\y)^{u_L}\z(\frac{\lambda_i}{\lambda_i+\sum_{l=1}^L\mu_i^l}\y).
\end{eqnarray}
  Therefore the first part of the theorem follows.

 Next we prove the second part of the theorem. One calculates from (\ref{up01}) and (\ref{up02}) that for $n\le -1,$
 $E_\omega(U_n)=M_{0}M_{-1}\cdots M_{n+1}.$ By Ward equation, we have that
 \begin{eqnarray*}
  E_\omega T_1 &=&E_\omega \xi_{0,1}+\sum_{i\le -1}E_\omega U_{i,1}E_\omega \xi_{i,k}+\sum_{i\le -1}E_\omega (U_{i,1}+...+U_{i,L})E_\omega \tilde{\xi}_{i+1,k}\\
   &=&\frac{1}{\lambda_0+\sum_{l=1}^L\mu_0^l}+\sum_{i\le-1}\frac{\mathbf{e}_1M_0\cdots M_{i+1}\mathbf{e}_1^T}{\lambda_i+\sum_{l=1}^L\mu_i^l}+\sum_{i\le-1}\frac{\mathbf{e}_1M_0\cdots M_{i+1}\mathbf{1}}{\lambda_{i+1}+\sum_{l=1}^L\mu_{i+1}^l}\\
   &=&\sum_{i\le0}\frac{1}{\lambda_i+\sum_{l=1}^L\mu_i^l}(\mathbf{e}_1M_0\cdots M_{i+1}\mathbf{e}_1^T+\mathbf{e}_1M_0\cdots M_{i}\mathbf1)\\
   &=&\sum_{i\le0}\frac{1}{\lambda_i}\mathbf{e}_1M_0\cdots M_{i+1}\mathbf1
    \end{eqnarray*}
 where the empty product equals to identity. \qed

\noindent{\it 3.2 LLN-Proof of Theorem \ref{lln} }

Once the quenched mean of $T_1$ has been calculated,  the proof of Theorem \ref{lln} follows basically as that in \cite{rit}.  Let $\eta_{n}=T_n-T_{n-1}$ for $n\ge 1.$  Then one follows from the stationarity of the environment that  $(\eta_n)_{n\ge 1}$ is a stationary sequence of random variables. Let $\nu_n$ be the number of states to the left of $n$ the process $N_t$ has ever visited between $T_n$ and $T_{n+1}.$ We have that
\begin{equation*}\begin{split}
  \lim_{n\rto}P(\nu_n\ge n/L)&\le\lim_{n\rto}P(\theta^{T_n}N_t\le(L-1)N/L\text{ for some }t<\theta^{T_n}T_{n+1})\\
&=\lim_{n\rto} P(N_t\le -n/L \text{ for some }t<T_1)\\
&=0,
\end{split}
  \end{equation*}
where the second line follows from the stationarity of the environment and the last line follows since $P(T_1<\infty)=1.$
For Borel set $A,B$ we have that
\begin{equation*}\begin{split}
 \lim_{n\rto}&P(\eta_1\in A, \eta_n\in B)=\lim_{n\rto}P(\eta_1\in A, \eta_n\in B, \nu_n<n/L)\\
&=\lim_{n\rto}P(\eta_1\in A)P( \eta_n\in B, \nu_n<n/L)=P(\eta_1\in A)P( \eta_n\in B),
\end{split}
  \end{equation*}
where the second equality follows because $\{\eta_1\in A\}\in \sigma\{\omega_i:i\le 0\}$  whereas $\{\eta_n\in B,\nu_n\le n/L\}\in \sigma\{\omega_i:i\ge 0\}.$ Thus we have shown that under probability $P,$ $\eta_n,{n\ge 1}$ are stationary and mixing. Then an application of Birkhoff's ergodic theorem yields that $P$-a.s.,
\begin{equation}\label{llt}
  \lim_{n\rto}\frac{T_n}{n}=E(T_1)=\mathbb ES(\omega).
\end{equation}

For  $t>0,$ there is a unique integer-valued random number $n_t$ such that $T_{n_t}\le t<T_{n_t+1}.$ We have that
\begin{equation}\label{jb}\frac{n_t-\nu_{n_t}}{T_{n_t+1}}\le \frac{N_t}{t}\le \frac{n_t+1}{T_{n_t}}\end{equation}
Suppose that $\mathbb ES(\omega)=\infty.$ Then we have from (\ref{llt}) and (\ref{jb}) that $P$-a.s., $\limsup_{t\rto}N_t/t\le 0.$
If $\gamma_L=0,$ then $N_t$ is recurrent and $\lim_{t\rto}N_t/t=0.$ If $\gamma_L<0.$ Then $N_t$ is transient to the right and $\liminf_{t\rto}N_t/t\ge 0.$ We conclude that whenever $\mathbb ES(\omega)=\infty,$ $P$-a.s., $\lim_{t\rto}N_t/t=0.$ Part (b) of the theorem is proved.

To prove part (a) of the theorem, suppose that $\mathbb ES(\omega)<\infty$ and define for $n\ge0,$ $\overline T_n=\inf\{k:\chi_k=n\}.$ Similarly as (\ref{llt}), $P$-a.s., $\lim_{n\rto}\overline T_n/n$ exists and is finite.
Since in every step, $\chi_n$ jumps at most a distance $L$ to the left, $0\le\nu_n\le L(\overline T_{n+1}-\overline T_n).$
Then we have that $P$-a.s., \begin{equation}\label{nul}0\le \lim_{n\rto}\nu_n/n\le \lim L(\overline T_{n+1}-\overline T_n)/n=0.\end{equation}

 Taking (\ref{llt}), (\ref{jb}) and (\ref{nul}) together, we have that  $P$-a.s., \begin{equation*}\lim_{t\rto}\frac{N_t}{t}=\mathbb (ES(\omega))^{-1}.\end{equation*}
Thus part (a) of the theorem is proved.\qed

\noindent{\large{\bf \large Acknowledgements:}} The author would like to thank Professor Wenming Hong for his useful comments on the paper.

 % \begin{center}
%{\section*{Acknowledgements}}
%\end{center}

\end{document}